# PRELIMINARIES ON BEST PROXIMITY POINTS IN CYCLIC MULTIVALUED MAPPINGS

M. De la Sen

*Institute of Research and Development of Processes. University of Basque Country

Campus of Leioa (Bizkaia) - Aptdo. 644- Bilbao, 48080- Bilbao. SPAIN

email: *manuel.delasen@ehu.es*

**Abstract**: This paper investigates the fixed points and best proximity points of multivalued cyclic self-mappings in metric spaces under a generalized contractive condition involving Hausdorff distances.

## 1. Introduction

Important attention is being devoted along the last years to investigate fixed point theory for multivalued mappings concerning some relevant properties like, for instance, stability of the iterations, fixed points of contractive and nonexpansive self-mappings and the existence of either common or coupled fixed points of several multivalued mappings or operators. This paper investigates some properties of fixed point and best proximity point results for multivalued cyclic self- mappings under a general contractive-type condition based on the Hausdorff metric between subsets of a metric space and which includes a particular case the contractive condition for contractive single-valued self-mappings, including the problems related to cyclic self-mappings. This includes strict contractive cyclic self -mappings and Meir-Keeler type cyclic contractions. Through this paper, we consider a metric space $(X,d)$ and a multivalued 2-cyclic self-mapping $T: A \cup B \to A \cup B$ (simply referred to as a multivalued cyclic self-mapping in the sequel) where $A$ and $B$ are nonempty closed subsets of $X$, so that $T(A) \subseteq B$ and $T(B) \subseteq A$ and $D = dist(A,B) \geq 0$. Let us consider the subset of the set of real numbers $\mathbf{R}_{0+} = \mathbf{R}_+ \cup \{0\} = \{z \geq 0 : z \in \mathbf{R}\}$, $\mathbf{R}_+ = \{z > 0 : z \in \mathbf{R}\}$, let the symbols $"\vee"$ and $"\wedge"$ denote the logic disjunction ("or") and conjunction ("and"), and define the functions:

$$M : (A \cup B) \times (B \cup A) \times [0,1) \times \Delta \to \mathbf{R}_{0+}, \quad \varphi(A \cup B) \times (B \cup A) \times [0,1) \times \Delta \to (0,1] \qquad (2.1)$$

as follows:

$$M(Tx,Ty,K,\alpha,\beta) = max[M_1(Tx,Ty,K), M_2(Tx,Ty,\alpha,\beta)]$$
$$= max\left[ K \, max\{d(x,y), d(x,Tx), d(y,Ty), 1/2(d(x,Ty)+d(y,Tx))\}, \alpha d(x,Tx) + \beta d(y,Ty) \right] \qquad (2.2)$$

$$\varphi(x,y,K,\alpha,\beta) = \begin{cases} 1 & \text{if } ([(\alpha,\beta) \in \Delta_1] \vee [(\alpha,\beta) \in \Delta_2]) \wedge (M_2(Tx,Ty,\alpha,\beta) > M_1(Tx,Ty,K)) \\ 1-\beta & \text{if } ([(\alpha,\beta) \in \Delta_3]) \wedge (M_2(Tx,Ty,\alpha,\beta) > M_1(Tx,Ty,K)) \\ \dfrac{1-\beta}{1-\beta+\alpha} & \text{if } ((\alpha,\beta) \in \Delta_4) \wedge (M_2(Tx,Ty,\alpha,\beta) > M_1(Tx,Ty,K)) \\ 1 & \text{if } (0 \leq K < 1/2) \wedge (M_2(Tx,Ty,\alpha,\beta) \leq M_1(Tx,Ty,K)) \\ 1-K & \text{if } (1/2 \leq K < 1) \wedge (M_2(Tx,Ty,\alpha,\beta) \leq M_1(Tx,Ty,K)) \end{cases} \qquad (2.3)$$

; $\forall x,y \in (A \cup B) \times (B \cup A)$ for some real constants $K$, $(\alpha,\beta) \in \Delta \subseteq [0,1) \times [0,1)$, where

$$\Delta = \{(\alpha,\beta) : \alpha \geq 0, \beta \geq 0, \alpha+\beta < 1\}$$

$$\Delta_1 = \{(\alpha,\beta) \in \Delta : \alpha \leq \beta, \alpha(1+\alpha) + \beta < 1\}; \quad \Delta_2 = \{(\alpha,\beta) \in \Delta : \alpha \geq \beta, \beta(1+\beta) + \alpha < 1\}$$

$$\Delta_3 = \left\{(\alpha,\beta) \in \Delta : \alpha \in (0, 1/2), \frac{1-\alpha}{2} > \beta \geq 1 - \alpha(1+\alpha)\right\}$$

$$\Delta_4 = \left\{(\alpha,\beta) \in \Delta : \alpha \in [0, 1/2], \frac{1-\alpha}{2} \leq \beta < 1 - \alpha(1+\alpha), \alpha(1+\alpha) + \beta(2-\beta) < 1\right\}$$

(2.4)

Note that $\varphi(A \cup B) \times (B \cup A) \times [0,1) \times \Delta \to (0,1]$ is non-increasing since all its partial derivatives with respect to $K, \alpha, \beta$ exist and are non-positive; $\forall x, y \in (A \cup B) \times (B \cup A)$ and note also that the subsets $\Delta_i \subset \Delta$; $i = 1, 2, 3, 4$ are disjoint by construction. Let $CL(X)$ be the family of all nonempty and closed subsets of the vector space $X$. If $A, B \in CL(X)$ then we can define $(CL(X), H)$ being the generalized hyperspace of $(X, d)$ equipped with the Hausdorff metric $H : CL(X) \to \mathbf{R}_{0+}$:

$$H(A, B) = \max\left\{\sup_{x \in A} d(x, B), \sup_{y \in B} d(y, A)\right\} \quad (2.5)$$

The distance between $A$ and $B$ is

$$D = dist(A, B) = \inf_{x \in A, y \in B} d(x, y) = \inf_{x \in A} d(x, B) = \inf_{y \in B} d(y, A) \quad (2.6)$$

## 2. Main results

The main result generalizes the contractive condition of multivalued self-mappings of [3] based on Hausdorff´s generalized metric. It includes, in particular, a contractive condition for multi-valued self-mappings. The main result follows:

**Theorem 2.1**: Let $(X, d)$ be a complete metric space and let $T : A \cup B \to A \cup B$ be a, in general, multivalued cyclic self-mapping, where $A, B \subset X$ are nonempty, closed and subject to the contractive constraint through $T : A \cup B \to A \cup B$:

$$[\varphi(x, y, K, \alpha, \beta) d(x, Tx) \leq d(x, y)] \Rightarrow [H(Tx, Ty) \leq M(Tx, Ty, K, \alpha, \beta) + \omega D] \quad (2.7)$$

for some $\omega \in \mathbf{R}_+$, $K \in [0, 1)$ and $(\alpha, \beta) \in \Delta$; $\forall (x, y) \in (A \cup B) \times (B \cup A)$. Assume also that

$$K_1 = \max\left(K, \frac{\beta}{1-\alpha}, \frac{\alpha}{1-\beta}\right) \in [0, 1), \quad K_2 = \max\left(\frac{1}{1-\alpha}, \frac{1}{1-\beta}\right) \in [1, \infty) \quad (2.8)$$

Then, the following properties hold:

**(i)** There is sequence $\{x_n\}$ in $A \cup B$ satisfying $x_{i+1} \in Tx_i$, $i \in \mathbf{N}$ such that

$$D \leq d(x_{n+1}, x_n) < \infty \; ; \; D \leq \limsup_{n \to \infty} d(x_{n+1}, x_n) \leq \frac{\omega K_2 D}{1 - K_1}$$

If $A$ and $B$ are bounded sets which intersect then $\sum_{n=1}^{\infty} d(x_{n+1}, x_n) < \infty$ and $\{x_n\}$ is a Cauchy sequence, then bounded of limit in $A \cap B$, with $x_{n+1} \in Tx_n$; $n \in N$ for any given $x_1 \in A \cup B$. If $A$ and $B$ are not bounded then the above property still holds if $d(x_1, x_2) < \infty$.

$\exists \lim_{n \to \infty} d(x_{n+1}, x_n) = D$ if $\omega = \dfrac{1 - K_1}{K_2}$ for any given $x_1 \in A \cup B$ with the sequence $\{x_n\}$ being constructed such that $x_{n+1} \in Tx_n$. If $x_1 \in A$ then the sequence of sets $\{Tx_{2n+1}\} \subset T(B) \subseteq A$ converges asymptotically to a subset $\{z_A\}$ of best proximity points in $A$ and the sequence of sets $\{Tx_{2n}\} \subset T(A) \subseteq B$ converges asymptotically to a subset $\{z_B\} \subset T\{z_A\} \subset T(A) \subset B$ of best proximity points in $B$ with $\{z_A\} \subset T\{z_B\} \subset T(B) \subset A$.

If $D = 0$, i.e. if $A \cap B \neq \emptyset$ then $\lim_{n \to \infty} \sup_{m > n} d(x_m, x_n) = \lim_{n \to \infty} d(x_{n+1}, x_n) = 0$, and any sequence $\{x_n\}$ being iteratively generated as $x_{n+1} \in Tx_n$, for any $x_1 = x \in A \cup B$, is a Cauchy sequence which converges to a fixed point $z \in Tz \cap (A \cap B)$ of $T : A \cup B \to A \cup B$. Assume that $D = 0$, i.e. $A \cap B \neq \emptyset$, that $A$ and $B$ are, in addition, convex and that $z_i \in Tz_i$; $i \in \overline{N} = \{1, 2, ...., N\}$ are fixed points of $T : A \cup B \to A \cup B$. Then, $z_i = z_j \subset A \cap B$ and $Tz_i \equiv Tz_j \subset A \cap B$; $\forall i, j (\neq i) \in \overline{N} = \{1, 2, ...., N\}$, that is, the image sets of any fixed points are identical.

(ii) Consider a uniformly convex Banach space $(X, \|\;\|)$, so that $(X, d)$ is a metric space for the norm - induced metric $d : X \times X \to R_{0+}$, and let $A$ and $B$ be nonempty, disjoint, convex and closed subsets of $X$ with $T : A \cup B \to A \cup B$ satisfying the contractive condition (2.7)-(2.8) with $\omega = \dfrac{1 - K_1}{K_2}$. Then, a sequence $\{x_{2n}\}$ built such that $x_{2n} \in Tx_{2n-1}$ with $x_{2n-1} \in Tx_{2n-2}$ is a Cauchy sequence in $A$ if $x_1 \in A$ and a Cauchy sequence in $B$ if $x_1 \in B$ so that $\lim_{n \to \infty} d(x_{2n+3}, x_{2n+1}) = \lim_{n \to \infty} d(x_{2n+2}, x_{2n}) = 0$; $\forall x_1 \in A \cup B$, and $\lim_{n \to \infty} d(x_{2n+2}, x_{2n+1}) = \lim_{n \to \infty} d(x_{2n+1}, x_{2n}) = D$; $\forall x_1 \in A \cup B$. If $x_1 \in A$ and $x_2 \in Tx_1 \subset T(A) \subset B$ then the sequences of sets $\{T^{2n} x_1\} \equiv \{T(T^{2n-1} x_1)\}$ and $\{T^{2n+1} x_1\}$ converge to best proximity points $z_A \in Tz_B$ and $z_B \in Tz_A$ in $A$ and $B$, respectively.

*Proof:* a) Take with no loss in generality $x \in A$ and $y \in Tx$ being arbitrary and assume that $\varphi(x, y, K, \alpha, \beta) d(x, Tx) \leq d(x, Tx) \leq d(x, y)$ since $\varphi(x, y, K, \alpha, \beta) \in (0, 1]$. Assume that $M_2(Tx, Ty, \alpha, \beta) > M_1(Tx, Ty, K)$ then,

$$d(y, Ty) \leq M_2(Tx, Ty) + \omega D = \alpha d(x, Tx) + \beta d(y, Ty) + \omega D \leq M_1(Tx, Ty) + \omega D$$
$$= K \max\left\{d(x, y), d(y, Ty), \dfrac{d(x, Ty)}{2}\right\} + \omega D = K \max\left\{d(x, y), \dfrac{d(x, Ty)}{2}\right\} + \omega D \qquad (2.9)$$

since $K \in [0,1)$ what also implies that $d(y, Ty) \leq \dfrac{1}{1-\beta}(d(x, Tx) + \omega D)$. Then,

$$d(y, Ty) \leq min\left(\dfrac{\alpha}{1-\beta} d(x, Tx), K\, max\left\{d(x, y), \dfrac{d(x, Ty)}{2}\right\}\right) + \dfrac{\omega D}{1-\beta}$$

$$\leq min\left(\dfrac{\alpha}{1-\beta} d(x, Tx),\, K\, d(x, y)\right) \leq max\left(K, \dfrac{\alpha}{1-\beta}\right) d(x, y) + \dfrac{\omega D}{1-\beta} \qquad (2.10)$$

since $d(x, Tx) \leq d(x, y)$ ; $\forall y \in Tx$.

b) Now, assume that $M_2(Tx, Ty, \alpha, \beta) \leq M_1(Tx, Ty, K)$ so that

$$d(y, Ty) \leq M_1(Tx, Ty) = K\, max\left\{d(x, y), d(y, Ty), \dfrac{d(x, Ty)}{2}\right\} = K\, max\left\{d(x, y), \dfrac{d(x, Ty)}{2}\right\} + \omega D$$

$$\leq M_2(Tx, Ty) = \alpha d(x, Tx) + \beta d(y, Ty) + \omega D \qquad (2.11)$$

This implies also that $d(y, Ty) \leq \dfrac{\alpha}{1-\beta} d(x, Tx) + \dfrac{\omega D}{1-\beta}$ and again (2.10) holds. As a result,

$$d(y, Ty) \leq max\left(K, \dfrac{\alpha}{1-\beta}\right) d(x, y) + \dfrac{\omega D}{1-\beta} \;\; ; \;\; \forall x \in A, \forall y \in B$$

by interchanging the roles of $A$ and $B$, one also gets by proceeding in a similar way:

$$d(x, Tx) \leq max\left(K, \dfrac{\beta}{1-\alpha}\right) d(x, y) + \dfrac{\omega D}{1-\alpha} \;\; ; \;\; \forall x \in A, \forall y \in B$$

Thus,

$$d(x, Tx) \leq max\left(K, \dfrac{\beta}{1-\alpha}, \dfrac{\alpha}{1-\beta}\right) d(x, y) + max\left(\dfrac{1}{1-\alpha}, \dfrac{1}{1-\beta}\right) \omega D = K_1 d(x, y) + K_2 \omega D \qquad (2.12)$$

; $\forall (x, y) \in (A \cup B) \times [T(A \cup B)]$, where $K_1 = max\left(K, \dfrac{\beta}{1-\alpha}, \dfrac{\alpha}{1-\beta}\right) \in [0,1)$ and

$K_2 = max\left(\dfrac{1}{1-\alpha}, \dfrac{1}{1-\beta}\right)$. Note that, since $T: A \cup B \to A \cup B$ is cyclic, then $y, Tx \in B$ if $x \in A$ and conversely.

Now, construct a sequence $\{x_n\}$ in $A \cup B$ as follows $x = x_1 \in A$, $x_2 \in Tx_1 \subset T(A) \subset B$,....
$x_{2n} \in Tx_{2n-1} \subset T(A) \subset B$,..., $x_{2n+1} \in Tx_{2n} \subset T(B) \subset A$ which satisfies:

$$D \leq d(x_{n+1}, x_n) \leq K_1 d(x_n, x_{n-1}) + \omega K_2 D \leq K_1^{n-1} d(x_2, x_1) + \left(\sum_{i=0}^{n-2} K_1^i\right) K_2 \omega D \;\; ; \;\; n \in N$$

$$\leq K_1^{n-1} d(x_2, x_1) + \dfrac{1 - K_1^{n-1}}{1 - K_1} K_2 \omega D < \infty \;\; ; \;\; n \in N \qquad (2.13)$$

Then, $D \leq \underset{n \to \infty}{lim\, sup}\, d(x_{n+1}, x_n) \leq \dfrac{K_2 \omega D}{1 - K_1}$. On the other hand,

$$\sum_{n=1}^{j} d(x_{n+1}, x_n) \leq \sum_{n=1}^{j} \left(K_1^{n-1} d(x_2, x_1) + \dfrac{1 - K_1^{n-1}}{1 - K_1} K_2 \omega D\right)$$

$$\leq \frac{1}{1-K_1}\left[\left(1-K_1^j\right)d(x_2,x_1)+\left(\sum_{n=1}^{j}\left(1-K_1^{n-1}\right)\right)K_2\omega D\right]; \ n\in N \qquad (2.14)$$

so that

$$\sum_{n=1}^{\infty}d(x_{n+1},x_n)\leq \frac{1}{1-K_1}\left[d(x_2,x_1)+\left(\sum_{n=1}^{\infty}\left(1-K_1^{n-1}\right)\right)K_2\omega D\right]; \ n\in N \qquad (2.15)$$

and we conclude that $\{x_n\}$ is a Cauchy sequence if $D=0$ (i.e. if $A$ and $B$ intersect provided that they are bounded or simply if $d(x_2,x_1)<\infty$), since $\lim\limits_{n\to\infty}\sup\limits_{m>n}d(x_m,x_n)=\lim\limits_{n\to\infty}d(x_{n+1},x_n)=0$, which has a limit $z$ in $X$, since $(X,d)$ is complete, which is also in $A\cap B$ which is nonempty and closed since $A$ and $B$ are both nonempty and closed since $T(A)\subseteq B$ and $T(B)\subseteq A$. On the other hand, for any distance $D\geq 0$ between $A$ and $B$:

$$D\leq d(x_{2n+3},x_{2n+2})\leq K_1^{2n+1}d(x_2,x_1)+\left(\sum_{i=0}^{2n}K_1^i\right)\omega D\leq K_1^{2n+1}d(x_2,x_1)+\frac{1-K_1^{2n+1}}{1-K_1}K_2\omega D\ ;\ n\in N \quad (2.16)$$

$$D\leq d(x_{2n+2},x_{2n+1})\leq K_1^{2n}d(x_2,x_1)+\left(\sum_{i=0}^{2n-1}K_1^i\right)\omega D\leq K_1^{2n}d(x_2,x_1)+\frac{1-K_1^{2n}}{1-K_1}K_2\omega D\ ;\ n\in N \qquad (2.17)$$

$$0\leq d(x_{2n+3},x_{2n+1})\leq K_1 d(x_{2n+2},x_{2n})+\omega D\leq K_1^2 d(x_{2n+1},x_{2n-1})+(1+K_1)K_2\omega D$$

$$\leq K_1^3 d(x_{2n},x_{2n-2})+(1+K_1+K_1^2)\omega D\leq K_1^4 d(x_{2n-1},x_{2n-3})+(1+K_1+K_1^2+K_1^3)K_2\omega D$$

$$\leq K_1^{2n}d(x_3,x_1)+\left(\sum_{i=0}^{2n-1}K_1^i\right)\omega D\leq K_1^{2n}d(x_3,x_1)+\frac{1-K_1^{2n}}{1-K_1}K_2\omega D\ ;\ n\in N \qquad (2.18)$$

Note that the sequences $\{d(x_n,x_{n+1})\}$ and $\{d(x_n,x_{n+2})\}$ are bounded if $x_1$ and $x_2\in Tx_1$ are such that $d(x_1,x_2)<\infty$ what is always guaranteed if $A$ and $B$ are bounded. If $\omega=\frac{1-K_1}{K_2}$ then one gets from the above relations that

$$\exists \lim_{n\to\infty}d(x_{2n+3},x_{2n+2})=\lim_{n\to\infty}d(x_{2n+2},x_{2n+1})=D$$

where $x_{2n+1}\in Tx_{2n}\subset A$, $x_{2n+2}\in Tx_{2n+1}\subset B$ and $x_{2n+3}\in Tx_{2n+2}\subset A$. Thus, any sequences of sets $\{x_{2n+1}\}$ and $\{x_{2n}\}$ contain asymptotically best proximity points of $A$ and $B$, respectively, if $x_1\in A$ and, conversely, of $B$ and $A$ if $x_1\in B$. This follows by contradiction since, if not, for each $k\in N$, there exists some $\varepsilon=\varepsilon(k)\in R_+$, some subsequence $\{n_{kj}\}_{j\in N}$ of natural numbers with $n_{km}>n_{kj}>k$ for $m>j$, and some related subsequences of real numbers $\{x_{2n_{kj}+1}\}$ and $\{x_{2n_{kj}}\}$ such that $d(x_{2n_{kj}+2},x_{2n_{kj}+1})\geq D+\varepsilon$ so that $d(x_{2n_k+2},x_{2n_k+1})\to D$ as $n_k\to\infty$ is impossible.

Now, assume $D=0$ and consider separately the various cases in (2.3)-(2.4), by using the contractive condition (2.7), subject to (2.2), to prove that there exists $z\in Tz$ in $A\cap B$ to which all sequences converge by using $D=0\Rightarrow \lim\limits_{n\to\infty}\sup\limits_{m>n}d(x_m,x_n)=\lim\limits_{n\to\infty}d(x_{n+1},x_n)=0\Rightarrow \{x_n\}\to z\in A\cap B$ with $\{x_n\}$ being a Cauchy sequence since $(X,d)$ is complete and $A$ and $B$ are nonempty and closed.

*Case 1*: $\varphi(x,y,K,\alpha,\beta)=1, ([(\alpha,\beta)\in \Delta_1]\vee [(\alpha,\beta)\in \Delta_2])\wedge (M(Tz,Tz,K,\alpha,\beta)=M_2(Tz,Tz,\alpha,\beta)>M_1(Tz,Tz,K))$.

Then, $d(z,Tz)=\varphi(z,z,K,\alpha,\beta)d(z,Tz)\leq (\alpha+\beta)d(z,Tz)\leq (1-\alpha^2)d(z,Tz)$ if $(\alpha,\beta)\in \Delta_1$. Thus, the contradiction $d(z,Tz)<d(z,Tz)$ holds if $(\alpha,\beta)\in \Delta_1, \alpha\neq 0$ and $z\notin Tz$. Hence, $z\in Tz$ if $(\alpha,\beta)\in \Delta_1$ with $\alpha\neq 0$ since $Tz$ is closed. If $\alpha=0$ then $0\leq \beta<1$ so that $d(z,Tz)\leq \beta d(z,Tz)<d(z,Tz)$ if $z\notin Tz$. Hence, $z\in Tz$ if $\alpha=0$ and $(0,\beta)\in \Delta_1$. The proof that $z\in Tz$ if $(\alpha,\beta)\in \Delta_2$ is similar since $(\alpha,\beta)\in \Delta_2 \Leftrightarrow (\beta,\alpha)\in \Delta_1$ from the definitions of the sets $\Delta_1$ and $\Delta_2$, and the fact that distances have the symmetry property.

*Case 2*: $\varphi(z,z,K,\alpha,\beta)=1-\beta, ([(\alpha,\beta)\in \Delta_3])\wedge (M(Tz,Tz,K,\alpha,\beta)=M_2(Tz,Tz,\alpha,\beta)>M_1(Tz,Tz,K))$.

Then, $(1-\beta)d(z,Tz)=\varphi(x,y,K,\alpha,\beta)d(z,Tz)\leq (\alpha+\beta)d(z,Tz)$ if $(\alpha,\beta)\in \Delta_3$ which fails for $z\notin Tz$ if and only if $1>\alpha+2\beta$. But $\Delta_3=\left\{(\alpha,\beta)\in \Delta: \alpha\in (0,1/2), \frac{1-\alpha}{2}>\beta\geq 1-\alpha(1+\alpha)\right\}$ so that $1>\alpha+2\beta$ and then the inequality fails if $z\notin Tz$. Hence, $z\in Tz$ since $Tz$ is closed.

*Case 3*: $\varphi(z,z,K,\alpha,\beta)=\frac{1-\beta}{1-\beta+\alpha}, ((\alpha,\beta)\in \Delta_4)\wedge (M(Tz,Tz,K,\alpha,\beta)=M_2(Tz,Tz,\alpha,\beta)>M_1(Tz,Tz,K))$.

Then, $\frac{1-\beta}{1-\beta+\alpha}d(z,Tz)=\varphi(z,z,K,\alpha,\beta)d(z,Tz)\leq (\alpha+\beta)d(z,Tz)$ if $(\alpha,\beta)\in \Delta_4$ which fails for $z\notin Tz$ if and only if $\frac{1-\beta}{1-\beta+\alpha}>\alpha+\beta$, equivalently if and only if $1>\alpha(1+\alpha)+\beta(2-\beta)$. But $\Delta_4=\left\{(\alpha,\beta)\in \Delta: \alpha\in [0,1/2], \frac{1-\alpha}{2}\leq \beta<1-\alpha(1+\alpha), \alpha(1+\alpha)+\beta(2-\beta)<1\right\}$ so that the inequality fails if $z\notin Tz$. Hence, $z\in Tz$ since $Tz$ is closed.

*Case 4*: $\varphi(x,y,K,\alpha,\beta)=1, (0\leq K<1/2)\wedge (M_2(Tz,Tz,\alpha,\beta)\leq M(Tz,Tz,K,\alpha,\beta)=M_1(Tz,Tz,K))$.

Then, $d(z,Tz)=\varphi(z,z,K,\alpha,\beta)d(z,Tz)\leq K\max\left\{d(z,z),d(z,Tz),\frac{d(z,Tz)}{2}\right\}=Kd(z,Tz)<d(z,Tz)$

which is a contradiction for any $z\notin Tz$. Hence, $z\in Tz$ since $Tz$ is closed.

*Case 5*: $\varphi(x,y,K,\alpha,\beta)=1-K, (1/2\leq K<1)\wedge (M_2(Tz,Tz,\alpha,\beta)\leq M(Tz,Tz,K,\alpha,\beta)=M_1(Tz,Tz,K))$. Then,

$(1-K)d(z,Tz)=\varphi(z,z,K,\alpha,\beta)d(z,Tz)\leq d(z,Tz)$

$$\Rightarrow d(z,Tz)\leq K\max\left\{d(z,z),d(z,Tz),\frac{d(z,Tz)}{2}\right\}=Kd(z,Tz)<d(z,Tz)$$

which is a contradiction if $z\notin Tz$. Hence, $z\in Tz$ since $Tz$ is closed. A combined result of Cases 1-5 is that $D=0\Rightarrow \{x_n\}\to z\in Tz\cap (A\cap B)$ for any $x_1\in A\cup B$. Now, assume again that $A\cap B\neq \emptyset$ and that there are two distinct fixed points $z_x(\neq z_y \in Tz_y) \in Tz_x$ necessary located in $A \cap B$ to which the sequences $\{x_n\}$ and $\{y_n\}$ converge to $z \in A \cap B$ and $q(\neq z) \in A \cap B$, respectively, where $x_{n+1} \in Tx_n$, $y_{n+1} \in Ty_n$ for $n \in N$ where $x_1, y_1(\neq x_1) \in A \cup B$. Assume also that $Tz \neq Tq$. One gets from the contractive condition (2.7), subject to (2.2)-(2.4), that:

$$\max(d(z, Tq), d(q, Tz)) \leq \max\left(\sup_{x \in Tz} d(x, Tq), \sup_{y \in Tq} d(y, Tz)\right)$$

$$\leq K \max(d(z, q), 1/2(d(z, Tq) + d(q, Tz))) = K d(z, q) < d(z, q)$$

Thus, construct sequences $z_{n+1} \in Tz_n$, $q_{n+1} \in Tq_n$ with $z_1 = z$ and $q_1 = q$ such that $d(z, q_{n+1}) < d(z, q_n)$ $d(q, z_{n+1}) < d(q, z_n)$ for $n \in N$. Since $z, q \in A \cap B$ which is nonempty, closed and convex, there exists $n_0 = n_0(\varepsilon)$ for any given $\varepsilon \in R_+$ such that $x_n$ and $q_n$ are in $A \cap B$ for $n \geq n_0$. Then, $q_n \to \hat{z}(\in Tz)$ and $z_n \to \hat{q}(\in Tq)$ as $n \to \infty$ with $z \in Tz \cap A \cap B$ and $q \in Tq \cap A \cap B$. Hence, $Tz \equiv Tq$ in $A \cap B$ contradicting the hypothesis that such sets are distinct. Property (i) has been proven.

Property (ii) is proven by using, in addition, [Lemma 3.8, ref. 4], one gets:

$$\exists \lim_{n \to \infty} d(x_{2n+3}, x_{2n+2}) = \lim_{n \to \infty} d(x_{2n+2}, x_{2n+1}) = D \Rightarrow \lim_{n \to \infty} d(x_{2n+3}, x_{2n+1}) = 0$$

for any sequence $\{x_n\}$ with $x_1 \in A \cup B$ and $x_{n+1} \in Tx_n$ since $(X, d)$ is a uniformly convex Banach space, $A$ and $B$ are nonempty and disjoint closed subsets of $X$ and $A$ is convex. Note that Lemma 3.8 of ref. 4 and its given proof remain fully valid for multivalued cyclic self-maps since only metric properties were used in its proof. It turns out that $\{x_{2n+1}\}$ is a Cauchy sequence, then bounded, with a limit $z_A$ in $A$, which is also a best proximity point of $T : A \cup B \to A \cup B$ in $A$ since

$$\lim_{n \to \infty} d(x_{2n+2}, x_{2n+1}) = \lim_{n \to \infty} d(x_{2n+2}, z_A) = D \leq \lim_{n \to \infty} d(x_{2n+2}, x_{2n}) + \lim_{n \to \infty} d(x_{2n}, z_A) = \lim_{n \to \infty} d(x_{2n}, z_A) \leq D$$

and then $\{x_{2n}\}$ converges to some point in $z_B \in Tz_A \subset B$, which is also a unique best proximity point in $B$ (then $z_B \in Tz_A$ and $Tz_A \subset B$), since $(X, d)$ is a uniformly convex Banach space and $A$ and $B$ are nonempty closed and convex subsets of $X$. In the same way, $z_A \in Tz_B \subset A$. Also, $\{x_{2n}\}$ is bounded and $\{x_{2n+1}\}$ is bounded sequences since $\{x_{2n}\}$ is bounded and $D < \infty$. Also, if $x_1 \in B$ and $B$ is convex then the above result holds with $x_{2n+1} \in Tx_{2n} \subset B$, $x_{2n+2} \in Tx_{2n+1} \subset A$ and $x_{2n+3} \in Tx_{2n+2} \subset B$. Now, for $D > 0$, the reformulated five Cases 1-5 in the proof of Property (i) would lead to contradictions $D = d(z_A, z_B) < D \neq 0$ if $z_A \notin Tz_B$ or if $z_B \notin Tz_A$. From Proposition 3.2 of ref. 4, there exist $z_A \in Tz_B$ and $z_B \in Tz_A$ such that $D = d(z_A, Tz_A) = d(z_B, Tz_B)$ since $T : A \cup B \to A \cup B$ is cyclic satisfying the contractive condition (2.7)-(2.8), where $A$ and $B$ are nonempty and closed subsets of a complete metric space $(X, d)$, with convergent subsequences $\{x_{2n+1}\}$ and $\{x_{2n}\}$ in both $A$ and $B$, respectively, for any $x = x_1 \in A$ and in $B$ and $A$, respectively, for any given $x = y_1 \in B$. Assume some given sequence $\{x_{2n+1}\}$ in $A$, being generated from some given $x_1 \in A$ with $x_{2n+1} \in Tx_{2n}$, which converges to the best proximity point $z_A \in A \cap Tz_B$ in $A$ of $T: A \cup B \to A \cup B$. Assume also that there is some sequence $\{y_{2n}\}$, distinct of $\{x_{2n}\}$, in $A$ generated from $y_1 (\neq x_1) \in A$ with $y_{2n+1} \in Ty_{2n}$ which converges to $z_{A1} \in A$ where $z_B \in B \cap Tz_A$ is a best proximity point in $B$ of $T: A \cup B \to A \cup B$. Consider the metric space $(X, d)$ obtained by using the norm-induced metric in the Banach space $(X, \|\ \|)$ so that both spaces can be mutually identified to each other. Since $d(x, y) \geq D$ for any $x \in A$ and $y \in B$, it follows that $D = d(z_B, z_A) = d(z_A, Tz_A) = d(z_B, Tz_B) < d(z_{A1}, Tz_A)$ if $z_{A1} \in A \cap \overline{Tz_B}$, where $z_A$ and $z_B$ are best proximity points of $T: A \cup B \to A \cup B$ in $A$ and $B$. Hence, $z_A, z_{A1} \in Tz_B$ and $z_B \in Tz_A$ and then any sequence converges to best proximity points. Hence Property (ii) has been proven. □

**ACKNOWLEDGMENT**S

The author is grateful to the Spanish Ministry of Education for its partial support of this work through Grant DPI2009-07197. He is also grateful to the Basque Government for its support through Grants IT378-10 and SAIOTEK S-PE09UN12.